%

%


\def\today{\ifcase\month\or January\or February\or
March\or April\or May\or June\or July\or August\or
September\or October\or November\or December\fi
\space\number\day, \number\year}




\def\dspace{\lineskip=2pt\baselineskip=18pt
\lineskiplimit=0pt}

\font \bbrm=cmbx10 at 12pt

\def\bigtype{\bbrm}

\hsize=13.5cm
\magnification=1200
\def\ce{\centerline}

\def\hb{\hfill\break}

\def\title #1{\null\bigskip\ce{\bigtype #1}
\bigskip}

\def\alp{\alpha}		
\def\bet{\beta}		
\def\gam{\gamma}		
\def\del{\delta}

\def\tet{\theta}		

\def\kap{\kappa}
\def\lam{\lambda}

\def\ome{\omega}		




    
\font\tenboldgreek=cmmib10
 \font\sevenboldgreek=cmmib10 at 7pt
\font\fiveboldgreek=cmmib10 at 7pt
\newfam\bgfam
\textfont\bgfam=\tenboldgreek
\scriptfont\bgfam=\sevenboldgreek
\scriptscriptfont\bgfam=\fiveboldgreek

\mathchardef\ggarrow="7010

\font\tengerman=eufm10 \font\sevengerman=eufm7
\font\fivegerman=eufm5
\font\tendouble=msym10 \font\sevendouble=msym7
\font\fivedouble=msym5

\textfont4=\tengerman \scriptfont4=\sevengerman
\scriptscriptfont4=\fivegerman
\newfam\dbfam
\textfont\dbfam=\tendouble \scriptfont\dbfam=
\sevendouble
\scriptscriptfont\dbfam=\fivedouble

\mathchardef\ng="702D
\mathchardef\dbA="7041
\mathchardef\sm="7072
\mathchardef\nvdash="7030
\mathchardef\nldash="7031
\mathchardef\lne="7008
\mathchardef\sneq="7024
\mathchardef\spneq="7025
\mathchardef\sne="7028
\mathchardef\spne="7029
\mathchardef\ltms="706E
\mathchardef\tmsl="706F

\mathchardef\dbA="7041


\mathchardef\dbA="7041 
\mathchardef\dbB="7042 
\mathchardef\dbC="7043 
\mathchardef\dbD="7044 
\mathchardef\dbE="7045 
\mathchardef\dbF="7046 
\mathchardef\dbG="7047 
\mathchardef\dbH="7048 
\mathchardef\dbI="7049 
\mathchardef\dbJ="704A 
\mathchardef\dbK="704B 
\mathchardef\dbL="704C 
\mathchardef\dbM="704D 
\mathchardef\dbN="704E 
\mathchardef\dbO="704F 
\mathchardef\dbP="7050 
\mathchardef\dbQ="7051 
\mathchardef\dbR="7052 
\mathchardef\dbS="7053 
\mathchardef\dbT="7054 
\mathchardef\dbU="7055 
\mathchardef\dbV="7056 
\mathchardef\dbW="7057 
\mathchardef\dbX="7058 
\mathchardef\dbY="7059 
\mathchardef\dbZ="705A

\def\sdp{\times \hskip -0.3em {\raise 0.3ex
\hbox{$\scriptscriptstyle |$}}} 


\def\cf{{\rm \,cf\,}}

\def\min{\mathop{\rm min}}










\def\ddownarrow{\big\downarrow \hskip-0.70em\raise
2pt\hbox {$\big\downarrow$}}
\def\longright #1#2 {\smash{\mathop{\hbox to
#1pt {\rightarrowfill}}\limits_{#2}}}
\def\sqr#1#2{{\vcenter{\hrule height.#2pt\hbox{\vrule
width.#2pt height#1pt \kern#1pt \vrule width.#2pt}
\hrule height.#2pt}}}

\def\buildrul#1\under#2{\mathrel{\mathop{\null#2}
\limits_{#1}}}

\def\boxit#1{\vbox{\hrule\hbox{\vrule\kern3pt
\vbox{\kern3pt#1 \kern3pt}\kern3pt\vrule}\hrule}}

\def\subheading#1{\medskip\goodbreak\noindent{\bf
#1.}\quad}

\def\sect#1{\goodbreak\bigskip\centerline{\bf#1}
\medskip}
\def\pr{\smallskip\noindent{\bf Proof:\quad}}
\def\onumber #1{\ooalign{\hfil\raise.07ex\hbox{
\hfill$\scriptstyle \,#1$\hfil}
\cr\cr{$\bigcirc$}}}
\def\onumber c{\ooalign{\hfil\raise.07ex\hbox
{\hfill$\scriptstyle \,c$\hfil}
\cr\cr{$\bigcirc$}}}
\def\alpcirc {\ooalign{\hfil\raise.07ex
\hbox{\hfill$\scriptstyle\alp\;$\hfill}\cr\cr
{$\bigcirc$}}}

\def\longmapright #1#2 {\smash{\mathop{\hbox to
#1pt {\rightarrowfill}}\limits^{#2}}}
\def\longmapleft #1 #2 {\smash{\mathop{\hbox to
#1 pt {\leftarrowfill}}\limits^{#2}}}

\def\references#1{\goodbreak\bigskip\par\centerline
{\bf References}\medskip\parindent=#1pt}
\def\ref#1{\par\smallskip\hang\indent\llap{\hbox
to \parindent{#1\hfil\enspace}}\ignorespaces}

\def\back{{\raise 2.5pt\hbox{$\,\scriptscriptstyle
\backslash\,$}}}
\def\bks{{\backslash}}
\def\part{\partial}
\def\lwr #1{\lower 5pt\hbox{$#1$}\hskip -3pt}
\def\rse #1{\hskip -3pt\raise 5pt\hbox{$#1$}}
\def\lwrs #1{\lower 4pt\hbox{$\scriptstyle #1$}
\hskip -2pt}
\def\rses #1{\hskip -2pt\raise 3pt\hbox
{$\scriptstyle #1$}}

\def\<#1{\left\langle{#1}\right\rangle}

\def\subinbn{{\subset\hskip-8pt\raise 0.95pt
\hbox{$\scriptscriptstyle\subset$}}}

\def\llvdash{\mathop{\|\hskip-2pt
\raise 3pt\hbox{\vrule height 0.25pt width 1.5cm}}}

\def\lvdash{\mathop{|\hskip-2pt \raise 3pt\hbox
{\vrule height 0.25pt width 1.5cm}}}

\def\fakebold#1{\leavevmode\setbox0=\hbox{#1}%
  \kern-.025em\copy0 \kern-\wd0
  \kern .025em\copy0 \kern-\wd0
  \kern-.025em\raise.0333em\box0 }

\font\msxmten=msxm10
\font\msxmseven=msxm7
\font\msxmfive=msxm5
\newfam\myfam
\textfont\myfam=\msxmten
\scriptfont\myfam=\msxmseven
\scriptscriptfont\myfam=\msxmfive
\mathchardef\rhookupone="7016
\mathchardef\ldh="700D
\mathchardef\leg="7053
\mathchardef\ANG="705E
\mathchardef\lcu="7070
\mathchardef\rcu="7071
\mathchardef\leseq="7035
\mathchardef\qeeg="703D
\mathchardef\qeel="7036
\mathchardef\blackbox="7004
\mathchardef\bbx="7003
\mathchardef\simsucc="7025

\def\rhookup{{\fam=\myfam \rhookupone}}

\def\bigsquare{{\fam=\myfam\bbx}}

\font\tencaps=cmcsc10
\def\smallcaps{\tencaps}

\def\author#1{\bigskip\ce{\smallcaps #1}\medskip}

\def\upddots{\mathinner{\mkern
1mu\raise 1pt \hbox{.}\mkern 2mu \mkern
2mu \raise 4pt\hbox{.}\mkern 1mu \raise 7pt\vbox
{\kern 7 pt\hbox{.}}} }

\def\varchi{\ooalign{{\raise
1.385pt\hbox{$\chi$}}\crcr\hbox{--}\crcr}}

\def\trianarrow{{\raise 2pt\hbox to 0.50cm
{\hrulefill}\triangleright}}

\null
\overfullrule=0pt
\sect{LESS SATURATED IDEALS}
\bigskip
\ce{by}
$$\vbox{\halign{\tabskip2em\hfil#\hfil&\hfil#\hfil\cr
\bf Moti Gitik&\bf Saharon Shelah
\footnote{${}^{1}$}{Partially supported by the Basic Research Fund,Israel Academy of Sciences.
Publication number 577.}\cr  
School of Mathematical Sciences&Hebrew University of
Jerusalem\cr
Sackler Faculty of Exact Sciences&Department of
Mathematics\cr
Tel Aviv University&Givat Ram, Jerusalem\cr
Ramat Aviv  69978 Israel&\cr}}$$

\dspace
\vskip1truecm
{\narrower\medskip
\subheading{Abstract}
We are proving the following:
\item{(1)} If $\kap$ is a weakly inaccessible then
$NS_\kap$  is not $\kap^+$-saturated.
\smallskip
\item{(2)} If $\kap$  is a weakly inaccessible and
$\tet <\kap$  is regular then $NS^\tet_\kap$  is
not $\kap^+$-saturated.
\smallskip
\item{(3)} If $\kap$  is singular then
$NS^{cf\kap}_{\kap^+}$ is not $\kap^{++}$-saturated. 

Combining this with previous results of Shelah, one
obtains the following:     
\item{(A)} If $\kap >\aleph_1$  then $NS_\kap$
is not $\kap^+$-saturated.
\item{(B)} If $\tet^+<\kap$  then $NS^\tet_\kap$
is not $\kap^+$-saturated.\medskip} 
\vskip1truecm
\def\drc{\diamond'_{\rm club}}

\sect{0.~~Introduction}

\footnote{${}^{}$}{Partialy supported by the Basic Research Fund,Israel
Academy of Sciences.Publication number 577.}
By a classical result of R. Solovay [So]
every stationary subset $S$  of a regular cardinal $\kap$
can be split into $\kap$ disjoint stationary subsets. 
Or in other terminology $NS_\kap\rhookup S$ the
nonstationary ideal over $\kap$ restricted to a
stationary subset $S$  of $\kap$ is not $\kap$-saturated.

A natural question is if it is possible to
replace $\kap$-saturatedness by $\kap^+$-saturatedness.
I.e. instead of completely disjoint $\kap$
stationary subsets to have $\kap^+$ with
pairwise nonstationary intersection.  K.~Namba
[Na] proved that $NS_\kap$  cannot be $\kap^+$-saturated
over a measurable $\kap$.  J.~Baumgartner, A.~Taylor and
S.~Wagon [Ba-Ta-Wa] improved this to a greatly Mahlo. 
S.~Shelah [Sh1,7] showed that
$NS_{\kap^+}\rhookup S$  is not $\kap^+$-saturated
if for some regular $\tet <\kap,\ \tet\not=
cf\kap\  S\cap \{\alp <\kap^+|cf \alp =\tet\}$
is stationary.  Actually a much more general
result is proved there.  For $\kap=\aleph_1$, it
is possible to have $NS_\kap\  \kap^+$-saturated
as was shown by J.~Steel and R.~Van Wesep [St-V]
from $AD_R$, by H.~Woodin [Wo1] from $AD$,  by
M.~Foreman, M.~Magidor and S.~Shelah [Fo-Ma-Sh]
from a supercompact and by S.~Shelah [Sh6] from
a Woodin cardinal.   
Also $NS_{\kap^+}\rhookup S$  can be
$\kap^{++}$-saturated for a regular $\kap$  and
$S\subseteq \{\alp <\kap^+\mid\cf\alp =\kap\}$
by T.~Jech and H.~Woodin [Je-Wo] building on K.
Kunen [Ku] construction of saturated ideal.  If
$\kap$  is a singular then by M.~Foreman [Fo]
and T.~Jech, H.~Woodin [Je-Wo] $NS_{\kap^+}\rhookup S$  
can be $\kap^+$-saturated for some $S\subseteq \{\alp
<\kap^+\mid cf\alp =cf\kap\}$.  

If $\kap$ is an inaccessible then T.~Jech and H.~Woodin
[Je-Wo] showed that $NS_\kap\rhookup$ Regulars can be
$\kap^+$-saturated and in [Gi2] was shown a
consistency of $NS_\kap \rhookup S$
$\kap^+$-saturated for $S\subseteq \kap$  such
that $S\cap\{\alp <\kap\mid cf\alp =\tet\}$  is
stationary for every regular $\tet <\kap$. 

The purpose of the present paper is to show that
$NS_\kap$  and $NS_\kap^\tet=NS_\kap\rhookup\{\alp
\mid cf\alp =\tet\}$ cannot be $\kap^+$-saturated
where $\kap$ is a weakly inaccessible and $\tet
<\kap$.  And for singular $\kap$, $NS_{\kap^+}^{cf\kap}$
cannot be $\kap^{++}$-saturated.     

Combining this with the previous result of Shelah [Sh1,7]
one gets the following:

\proclaim Theorem.
\item{(1)} If $\kap >\aleph_1$  then $NS_\kap$
is not $\kap^+$-saturated
\item{(2)} If $\tet^+<\kap$  then $NS^\tet_\kap$
is not $\kap^+$-saturated.

The proof is based on a certain combinatorical
principle which follows from saturatedness.  It
was considered independently by both authors.
Actually they came to it from different directions.
Thus Gitik was dealing with indiscernables of the
Mitchell Covering Lemma and showed it first
under $\neg\exists \kap o(\kap)=\kap^{++}$ in
[Gi3] and then removed this assumption in [Gi4].
Shelah [Sh2,3,8] went through his club
guessing machinery and eventually obtained a
better result not appealing to some GCH
assumptions used in [Gi4].  He also realized
that this leads to inconsistency of saturatedness
of $NS_\kap$  over small inaccessibles.  A
minor change in this argument by Gitik led to
the full result for inaccessibles.Shelah saw that it gives also the
nonsaturatedness of $NS_{\kap^+}^{cf \kap}$ for a singular $\kap$.	 
\vfill\eject

\sect{1.~~Main Result}

Let $\kap$  be a regular uncountable cardinal.
We denote by $NS_\kap$  the nonstationary ideal
over $\kap$.  For a set $S\subseteq\kap$  let
$NS_\kap\rhookup S$  denote the collection of
all subsets of $\kap$  having a nonstationary
intersection with $S$.  For a regular $\tet
<\kap$ let 
$$S^\tet_\kap =\{\alp <\kap \mid cf\alp =\tet\}\ .$$
We denote $NS_\kap\rhookup S_\kap^\tet$  simply
by $NS^\tet_\kap$.

An ideal $I$ over $\kap$  is called $\lam$-saturated
if there is no sequence $\langle A_\alp \mid\alp
<\lam\rangle$ so that
\item{(a)} $A_\alp\subseteq\kap$
\smallskip
\item{(b)} $A_\alp\notin I$
\smallskip
\item{(c)} $A_\alp\cap A_\bet\in I$ whenever
$\alp\not=\bet$.
\item{}

\smallskip
Our main objective will be the $\kap^+$-saturatedness
of $NS^\tet_\kap$.So further by saturated we will mean
$\kap^+$ saturated.

 We will prove the following:

\proclaim Theorem 1.  
\item{(1)} If $\kap$  is a weakly inaccessible
then $NS_\kap$  cannot be saturated.
\item{(2)} If $\kap$  is a weakly inaccessible
and $\tet <\kap$  is a regular then
$NS^\tet_\kap$  cannot be saturated.
\item{(3)} If $\kap$ is a singular cardinal then
$NS^{\cf\kap}_{\kap^+}$  cannot be saturated.

The proof will consist of two parts.  The first
will introduce a combinatorical principle and
show its inconsistency.  The second will be to
show that the saturatedness of $NS_\kap^\tet$
implies the principle.

\subheading{Definition 2} For regular cardinals
$\kap,\tet$  with $\kap >\aleph_2$  and
$\kap >\tet^+$  let $\diamond^*_{{\rm club}}(\kap,\tet)$
denote the following:\hb 
there exists a sequence
$\langle S_\alp\mid\alp\in S^\tet_\kap\rangle$
so that for every $\alp\in S^\tet_\kap$
\item{(1)} $S_\alp\subseteq \alp$
\smallskip
\item{(2)} $\sup S_\alp =\alp$
\smallskip
\item{(3)} $|S_\alp|=\tet$
\smallskip
\item{(4)} for every $\bet\in S_\alp$ $cf\bet >\tet$
and if $\tet =\aleph_0$ then $cf\bet >\aleph_1$.
\smallskip
\item{(5)} for every club $C\subseteq\kap$ the set
$$\{\alp\in S^\tet_\kap\mid\exists\bet <\alp\
C\supseteq S_\alp\bks\bet\}$$
contains a club intersected with $S^\tet_\kap$.

This is a strengthening of Shelah's club guessing
principles see [Sh2,3].  It turns out that it is
too strong.  Namely:

\proclaim Lemma 3.
$ZFC\vdash^\neg\diamond^*_{\rm club}(\kap,\tet)$.

\pr  Suppose otherwise.  Let $\langle
S_\alp\mid\alp\in S^\tet_\kap\rangle$ be a
sequence witnessing $\diamond^*_{{\rm
club}}(\kap,\tet)$.

Let us split the proof into two cases according
to $\tet >\aleph_0$  or $\tet =\aleph_0$.

\subheading{Case 1}
$\tet >\aleph_0$. 

We define by induction a sequence $\langle
E_n\mid n<\ome\rangle$  of clubs of $\kap$.  Let
$E_0=\kap$.  If $E_n$  is defined then consider
$E'_n$  the set of limit points of $E_n$.  There
exists a club $C\subseteq\kap$ such that for every
$\alp\in C\cap S_\kap^\tet$  $E'_n$ contains a
final segment of $S_\alp$.  Set $E_{n+1}=E'_n\cap C$.

Finally let $E=\cap_{n<\ome} E_n$.  Then $E$  is
a club. Let $\del=\min \big(E\cap S^\tet_\kap\big)$.
For every $n<\ome$  $E'_n$ contains a final
segment of $S_\del$.  Since $\del\in E_{n+1}$.
But $cf\del =\tet>\aleph_0$. So the final
segment of $S_\del$  is contained in $E$.  Pick
some $\bet\in E\cap S_\del$. By Definition 2(4),
$cf\bet >\tet$.  Since $E'_n\supseteq E_{n+1}$, $\bet\in
E'_n$  for every $n<\ome$.  So it is a limit
point of $E_n$.  Hence $E_n\cap\bet$  is a club
of $\bet$,  for every $n<\ome$.  Since $\cf \bet
>\aleph_0$, also $E\cap\bet$  is a club
of $\bet$.  But $\cf\bet >\tet$.  Hence there is
some $\gam\in E\cap\bet$  of cofinality $\tet$
which contradicts the minimality of $\del$.

\subheading{Case 2} $\tet =\aleph_0$ 

We define a decreasing sequence of clubs as
above but of the length $\aleph_1$.  Le $\langle
E_i\mid i<\ome_1\rangle$ be such a sequence and
let $E=\bigcap_{i<\ome_1} E_i$. As in Case 1 pick
$\del$  to be $\min\big(E\cap S^{\aleph_0}_\kap\big)$.
By Definition 2(3), $|S_\del| =\ome$.  Let
$\langle s_n\mid n<\ome\rangle$  be a cofinal in
$\del$  sequence of $S_\del$. Since for
every $i<\ome_1$  $E_{i+1}$  contains a final
segment of $S_\del$, there will be some
$n^*<\ome$  such for $\aleph_1$  $i$'s
$E_{i+1}\supseteq \{s_n\mid n\ge n^*\}$.
But the sequence $\langle E_i\mid i<\ome_1\rangle$
is a decreasing sequence.  So for every
$i<\ome_1$  $E_i\supseteq \{s_n\mid n\ge
n^*\}$.  Hence there is $\bet\in E\cap S_\del$.
Notice that by Definition 2(4), $cf\bet
>\aleph_1$.  Now we continue toward the
contradiction as in Case 1.\hfill$\bigsquare$ 

\proclaim Lemma 2.  Suppose that $\kap,\tet$
are regular cardinals and $\tet^+ <\kap$.       

If $NS^\tet_\kap$  is saturated then
$\diamond^*_{\rm club}(\kap, \tet)$  holds.

\subheading{Remark}  Both authors arrived to
this statement independently and from different
directions.  Gitik dealing with indiscernibles
of Mitchell's Covering Lemma showed this under
$\neg\big(o(\alp)=\alp^{++}\big)$  and $2^\tet
<\kap$  in [Gi3] and much later in [Gi4] realized that
$\neg\big(o(\alp)=\alp^{++}\big)$ is not needed.
Shelah went through his club guessing principles
[Sh2,3] and got a better result requiring only
$\tet^+ <\kap$  and not $2^\tet <\kap$.  The
proof below follow his lines.

\pr Let $\diamond'_{\rm club}(S)$  for $S$  a
stationary subset of $S^\tet_\kap$  be defined
as $\diamond^*_{\rm club}(\kap , \tet)$ only
with $S$  replacing $S^\tet_\kap$  and of
Definition 2(5) replaced by the following:\hb
${\rm (5)}^\prime$  for every club
$C\subseteq\kap$ the set
$$\{\alp\in S\mid\exists\bet <\alp\quad
C\supseteq S_\alp\bks\bet\}$$
is stationary. 

It was shown by Shelah [Sh2,3] that $\diamond'_{\rm
club}(S)$  and even a stronger principles are
true in ZFC.  For the benefit of the reader we
will present a proof of $\diamond'_{\rm club}(S)$
below.  But first let us use it to complete the proof
of the lemma.

\subheading{Claim 2.1}  For every stationary
$S\subseteq S^\tet_\kap$  there exists a
stationary $S^*\subseteq S$  such that for every
club $C\subseteq\kap$  the set
$$\{\alp\in S^*\mid\exists\bet <\alp\quad
C\supseteq S_\alp\bks\bet\}$$
contains a club intersected with $S^*$ where
$\langle S_\alp\mid\alp\in S\rangle$ is a sequence
witnessing $\diamond'_{\rm club}(S)$. 

\pr Suppose otherwise.  We define by induction
an almost (modulo $NS_\kap$) decreasing
sequence $\langle C_\alp\mid\alp <\kap^+\rangle$
of clubs of $\kap$  and almost disjoint sequence
$\langle A_\alp\mid\alp <\kap^+\rangle$ of
stationary subsets of $S$.  Denote for a club
$C$ by $N(C)$  the set of places where $C$  is
guessed, i.e.
$$\{\alp\in S\mid\exists\bet <\alp\quad
C\supseteq S_\alp\bks \bet\}\ .$$
Pick $C_0$  to be a club such that $S\bks N(C_0)$
is stationary.  Set $A_0=S\bks N(C_0)$. Suppose now that
$\langle C_\bet\mid\bet <\alp\rangle$  and
$\langle A_\bet\mid\bet <\alp\rangle$ are defined. 
We like to define $C_\alp$  and $A_\alp$. First take a
club $C$  which is almost contained in every $C_\bet$ 
for $\bet <\alp$.  Consider $N(C)$.  Obviously $N(C)$  is
almost contained in every $N(C_\bet)$  for $\bet <\alp$.
Let $C_\alp$  be a club subset of $C$ witnessing that
$N(C)$  is not good, i.e.  $N(C)\bks N(C_\alp)$  is
stationary.  Set $A_\alp =N(C)\bks N(C_\alp)$.  
This completes the inductive definition.

The existence of the sequence $\langle A_\alp\mid\alp
<\kap^+\rangle$  contradicts the
saturatedness of $NS_\kap^\tet$.\hb
\hfill$\bigsquare$ of the claim 

Now using saturatedness, it is easy to glue
together the sets $S^*$ given by Claim 2.1 and
to get
$\diamond^*_{\rm club}(\kap,\tet)$.  Just pick a maximal
pairwise almost disjoint collection of such
$S^*$'s.  Since $NS^\tet_\kap$  is saturated, it
consists of at most $\kap$  many sets.  Then make
them completely disjoint using the normality of
$NS^\tet_\kap$.  Finally, we put together the
sequences witnessing $\drc$ and obtain
$\diamond^*_{\rm club}(\kap,\tet)$.
This completes the proof of the
lemma.\hfill$\bigsquare$

\subheading{Remark}  We can replace the
$NS_\kap^\tet$ by any normal ideal over $\kap$
concentrating on $S^\tet_\kap$  in Lemma 2.  Only the
definition of $\diamond^*_{{\rm club}}(\kap,\tet)$ 
be changed in the obvious way.  The
same proof works.  We refer to [Gi4] for generalizations
for ideals preserving $2^\tet$ and to Dzamonja, Shelah
[Dz-Sh] for further generalizations in this direction. 

Let us now present a proof of Shelah's club guessing
principle, see [Sh2,3], [Br-Ma] for variations of it.   

\proclaim Proposition 3.  Suppose that $\kap$,
$\tet$  are regular and $\kap >\max (\tet^+,
\aleph_2)$.  Then $\drc(S)$  holds for every
stationary $S\subseteq S^\tet_\kap$.    

\pr Suppose for simplicity that
$\tet^+\ge\aleph_2$.  Otherwise we do the same
only with $\tet^+$ replaced by $\tet^{++}$.  

Suppose that for some stationary $S\subseteq
S^\tet_\kap$  $\drc (S)$  fails.
We shall define sequences $\langle C_i\mid
i<\tet^+\rangle$  of clubs of $\kap$, $\langle
T^i_\alp \mid i<\tet^+$,  $\alp\in S\rangle$
of trees and $\langle S^i_\alp\mid i <\tet^+,
\alp\in S\rangle$.  However, first for a club
$C$  and an ordinal $\alp\in S^\tet_\kap$  let us
define canonically a tree $T_\alp (C)$.  The
first level of $T_\alp (C)$ will consist of a
closed cofinal in $\alp$  sequence of order type
$\tet$  so that each nonlimit point of it has
cofinality $>\tet$.  We pick such a sequence to
be the least in some fixed well-ordering.   

Now let $\eta$  be a point from the first level.
We like to define the set of its immediate
successors, i.e, the set $Suc_{T_\alp(C)}(\eta)$.
If $\eta$  is a limit point then
$Suc_{T_\alp(C)}(\eta)=\emptyset$.  Otherwise, let
$\eta^*$  be the largest point of the first level below
$\eta$.  Consider $\eta'=\sup (C\cap\eta +1)$.  If
$\eta'=\eta$  or $\eta'\le\eta^*$  then set
$Suc_{T_\alp(C)}(\eta)= \emptyset$.  Suppose that
$\eta^*<\eta'<\eta$.  If $cf\eta'>\tet$  then
set $Suc_{T_\alp(C)}(\eta)=\{\eta'\}$.
Such $\eta'$  will be a leaf, i.e. $Suc_{T_\alp(C)}
(\eta')=\emptyset$.  If $cf\eta'>\tet$  then
as above, we pick the least closed cofinal in $\eta'$ 
sequence of order type of $\cf\eta'$  with nonlimit points
of cofinality $>\tet$ and the first element
above $\eta^*$.  The set of immediate
successors of $\eta$ will consist of this
sequence.  Using it we continue to define $T_\alp(C)$
above $\eta$ in the same fashion. 

Obviously such a defined tree $T_\alp(C)$  is well
founded and of cardinality $\le\tet$.

Now let us turn to the definition of the sequences.
Set $C_0=\kap$.  For $\alp\in S$  set $T^0_\alp
=T_\alp (C_0)$  and $S^0_\alp =$ the set of all
points of all the levels of $T^0_\alp$
which are in $C_0$ and has cofinality $>\tet$.  Clearly,
for all but nonstationary many $\alp$'s in $S$, $S_\alp$
is unbounded in $\alp$.

Now let $C_1$  be a club subset of $C_0$
witnessing the failure of $\langle
S^0_\alp\mid\alp\in S\rangle$  to be a $\drc(S)$
sequence.  We define $T^1_\alp$  and $S^1_\alp$
as above for $C_1$  replacing $C_0$.  Continue
by induction. At limit stages $i$ we take $C_i$
to be $\bigcap_{j<i} C_j$.   

Finally, let $D=\bigcap_{i<\tet^+}C_i$.  Then
$D$ is a club since $\kap >\tet^+$. Let $\alp\in
D\cap S$  be so that elements of $D$  of
cofinality $>\tet$  are unbounded in it.  Let us
show that the trees $T_\alp(C_i)$  $(i<\tet^+)$
must stabilize.  I.e. starting with some $i_\alp
<\tet^+\ T_\alp (C_i)=T_\alp(C_j)$. Otherwise, pick some
$\eta_0\in Lev_1(T^0_\alp)$  such that there is no
stabilization above it. 

Notice that the first level in all the trees is
the same.  Then $\eta_0$  cannot be a limit
point since otherwise it will have the empty set
of successors.  Let $\eta^*_0$  be the largest
point of the first level below $\eta_0$.  Since
there is no stabilization above $\eta_0$  in
the trees and $C_i$'s are decreasing, the
sequence $\langle \eta_{0i}\mid i<\tet^+\rangle$
will be a nonincreasing sequence of ordinals
inside the interval $(\eta_0^*,\eta_0]$  where
$\eta_{0i}=\sup\big(C_i\cap (\eta_0+1)\big)$.
Hence it is eventually constant.  So there is
some $\eta_1\in (\eta^*_0,\eta_1)$,
$cf\eta_1\le\tet$  such that $\eta_{0i}=\eta_1$
starting with some $i(1)<\tet^+$.  By the
definition of the trees the set of immediate
successors of $\eta_1$  will be the same in
every $T_\alp(C_i)$  $(i\ge i(1))$.  Pick
$\eta_2$  to be one of them with no
stabilization above it.  Deal with it as it was
done with $\eta_0$.  We will obtain $\eta_3$
and $\eta_4$.  Continue in the same fashion.
This process will produce an infinite decreasing
sequence of ordinals. 

Hence, there is $i_\alp <\tet^+$  such that
$T^i_\alp =T_\alp^{i_\alp}$  for every $i\ge
i_\alp$.  Then there are $i^*<\tet^+$  and a
stationary $S^*\subseteq S$ such that for every
$\alp\in S^*$ $i_\alp =i^*$  and $\sup
S_\alp^{i^*}=\alp$.  But this contradicts the
choice of $C_{i^*+1}$.  Contradiction.\hfill$\bigsquare$

\sect{2.~~Some Open Problems On the Nonstationary
Ideal} 

The following is probably the most interesting
left open problem on this subject:

\subheading{Problem 1}  Let $\kap$  be a regular
uncountable cardinal.  Can $NS^\kap_{\kap^+}$
be saturated?

By Shelah [Sh1] the answer is negative if we replace
$NS^\kap_{\kap^+}$ by $NS^\tet_{\kap^+}$  some
$\tet <\kap$.

Now for $\kap =\aleph_0$,  the models with
$NS_{\aleph_1}$  saturated were constructed by
J.~Steel and R.~Van Wesep [St-V], H.~Woodin [Wo]
using $AD_R$  and $AD$  then by M.~Foreman,
M.~Magidor and S. Shelah [Fo-Ma-Sh] from a
supercompact and by S.~Shelah [Sh6] from a Woodin
cardinal.  By J.~Steel [St] ``There is a saturated
ideal over $\aleph_1$  and a measurable" implies
an inner model with a Woodin cardinal.
However, the following basic question remains
open:

\subheading{Problem 2}  Is it consistent $GCH
+NS_{\aleph_1}$  is saturated? 

A very tight connection of $\neg CH$  with saturatedness
of $NS_{\aleph_1}$  was established by H.~~Woodin.
 He showed that ``$NS_{\aleph_1}$ is
saturated and there is a measurable" implies
$\neg CH$. 

If we relax our assumptions and
consider saturatedness of $NS_\kap\rhookup S$
for a stationary set $S$, then by T.~Jech and H.~Woodin 
for every $\kap$ 
$NS_{\kap^+}\rhookup S$  can be saturated 
.  Also $NS_\kap\rhookup$  Regulars
can be saturated for inaccessible $\kap$. An almost
huge is used for the first result and a measurable for the
second.  By [Gi2], $NS_\kap\rhookup S$ can be
saturated over inaccessible with $S$  having a
stationary intersection with every $S^\tet_\kap$
for $\tet$  regular less than $\kap$ there $o(\kap)=\kap$
is used for this and it is necessary. But in
order to obtain the same over the
first inaccessible, supercompacts are used.  

\subheading{Problem 3}  How strong is $NS_\kap\rhookup
S$  saturated for the first inaccessible $\kap$
and $S\subseteq S_\kap^\tet$  $\tet <\kap$
regular?

One can try to get simultaneously many $S$'s for
which $NS_\kap\rhookup S$  is saturated.  So the
following is natural:  we call $NS_\kap$
densely saturated if for every stationary
$S\subseteq \kap$  there is $S^*\subseteq S$
such that $NS_\kap\rhookup S^*$  is saturated. 

\subheading{Problem 4} Can $NS_\kap$ be densely
saturated over an inaccessible $\kap$?
Can $NS_{\kap^+}^{cf\kap}$ be densely saturated for any uncountable $\kap$?

Further weakenings are leading to notions of
presaturatedness and precipitousness introduced
by J.~Baumgartner-A.~Taylor [Ba-Ta] and
T.~Jech-K.~Prikry [Je-Pr] respectively. 

An ideal $I$  over $\kap$  is called
presaturated if the forcing with it preserves
$\kap^+$.  It is precipitous if the generic ultrapower
is well-founded.  Saturatedness implies presaturatedness
and presaturatedness implies precipitousness.   
It turned out that for every $\kap$  $NS_\kap$
can be precipitous.  With presaturatedness the
situation is less clear.  Namely, $NS_\kap$  can
be presaturated for $\kap =\aleph_1$  or $\kap$
an inaccessible but by [Sh1] these are the only
cases.
See Jech-Magidor-Mitchell-Prikry [Je-Ma-Mi-Pr], [Gi1],
[Gi2], [Gi3], [Gi4], M.~Foreman, M.~Magidor, S. Shelah
[Fo-Ma-Sh] for the consistency results. 
We do not know the following 

\subheading{Problem 5}  Can $NS^{\aleph_1}_{\aleph_2}$
be presaturated?

Another direction is to consider
$\diamond_\kap$. Obviously $\diamond_\kap\longrightarrow
(NS_\kap\  \hbox{is not saturated})$.  One
can try to prove $\diamond_\kap$  in $ZFC$+ instances of $GCH$, of course.
Since by Foreman and Woodin [Fo-Wo] $GCH$ can fail everywhere and so is
$\diamond_\kap$.
H.Woodin [Cu-Wo] showed a consistency of
$\neg\diamond_\kap$ over the
first strongly inaccessible.
But $GCH$ fails on a club of $\kap$.
He started from a supercompact.  R.~Jensen proved that
at least measurable is needed for
this.

\subheading{Problem 6} Is it consistent $\kap$ is an inaccessible,
$\neg\diamond_\kap$ and $GCH$ ?

By [Sh 5],it is possible to have a stationry subset $S$ of an
inaccessible $\kap$ with $\neg\diamond_S$ and $GCH$.

\subheading{Problem 7} Is it consistent $2^\kap$=$\kap^+$ and
$\neg\diamond_{S_{\kap^+}^{cf\kap}}$ for a singular $\kap$?
The same for
$\neg\diamond_{\kap^+}$ ?

It is known that large cardinals are needed for this,see [Sh5].
By Gregory[Gr] and [Sh4] $\diamond_{S_\lambda^\kap}$ holds for a 
cardinal $\lambda=2^\mu=\mu^+$ and regular $\kap<\mu$ provided
$\mu^\kap=\mu$ or $\mu$ is singular $\kap$ differs from $cf\mu$
and for every $\del<\mu$, $\del^\kap<\mu$.
This was continued in [Sh 9].
Thus it was shown that for $\kap$ above $\bet_\omega$ $\diamond_
{\kap^+}$ is equivalent to $2^\kap$=$\kap^+$.Also see
 [Sh5] and Dzamonja,Shelah [Dz-Sh2] for related results.
On the otherhand Shelah showed that 
it is consistent $GCH$ and $\neg\diamond_{S_{\aleph_2}^{\aleph_1}}$
see Kim-Steinhorn[Ki-St].

\references {68}
\smallskip
\ref{[Ba-Ta]} J. Baumgartner and A. Taylor,
Saturation properties of ideals in generic
extensions I, II Trans. Am. Math Soc. 270
(1982), 557-574 and 271 (1982), 587-609. 

\smallskip
\ref{[Ba-Ta-Wa]} J. Baumgartner, A. Taylor and
S. Wagon, On splitting stationary subset of large
cardinals, J. Sym Logic 42 (1977), 203-214.

\smallskip
\ref{[Br-Ma]} M. Bruke and M. Magidor, Shelah's
pcf theory and its applications, Ann. Pure and
Ap. Logic 50 (1990), 207-254. 

\smallskip
\ref{[Cu-Wo]} J. Cummings and H. Woodin,
Generalized Prikry Forcing.

\smallskip
\ref{[Dz-Sh]} M. Dzamonja and S. Shelah, Some
results on squares, outside guessing of clubs
and $I_{<f}[\lam]$.  

\smallskip
\ref{[Dz-Sh2]} M.Dzamonja and S.Shelah,
Saturated filters at successors of singulars,weak reflection.
545.

\smallskip
\ref{[Fo]} M. Foreman, More saturated ideals, in
Cabal Seminar 79-81, Lecture Notes in Math. 1019
Springer 1983, 1-27.

\smallskip
\ref{[Fo-Ma-Sh]} M. Foreman, M. Magidor and S.
Shelah, Martin Maximum, saturated ideals and
nonregular ultrafilters, I, Ann. of Math. 127
(1988), 1-47.  

\smallskip
\ref{[Fo-Wo]} M.Foreman and H.Woodin,GCH can fail everywhere,
Ann. of Math. 133 (1991) , 1-36.

\smallskip
\ref{[Gi1]} M. Gitik, The nonstationary ideal on
$\aleph_2$,  Israel Journal of Math. V.58, No.4
(1984), 257-288.

\smallskip
\ref{[Gi2]} M. Gitik, Changing cofinalities and
the nonstationary ideal, Israel Journal of Math.
V.56 (1986), 280-314.  

\smallskip
\ref{[Gi3]} M. Gitik, Some results on the
nonstationary ideal, Israel Journal of Math., to
appear.

\smallskip
\ref{[Gi4]} M. Gitik, Some results on nonstationary
ideal 2, submitted to Israel J. of Math.

\smallskip
\ref{[Gr]} J. Gregory, Higher Souslin trees and
GCH, J. Symb. Logic 41 (1976), 663-671.

\smallskip
\ref{[Je-Wo]} T. Jech and H. Woodin,
Saturation of the closed unbounded filter on the
set of regular cardinals, Trans. Am. Math. Soc.
292 (1), (1985), 345-356.  

\smallskip
\ref{[Je-Ma-Mi-Pr]} T. Jech, M. Magidor,W.Mitchell and K.
Prikry, Precipitous ideals, J. Symb. Logic 45
(1980), 1-8.

\smallskip
\ref{[Je-Pr]} T.Jech and K.Prikry,
Ideals of sets and the power set operation,
Bull. Am. Math. Soc. 82 (1976), 593-595.

\smallskip
\ref{[Ki-St]} Kim and C. Steinhorn, Is.J. of Math.

\smallskip
\ref{[Ku]} K. Kunen, Saturated ideals, J. of
Sym. Logic 43 (1978), 65-76.

\smallskip
\ref{[Na]} K. Namba, On closed unbounded ideal
of ordinal numbers, Comm. Univ. Sancti Pauli
22 (1974), 33-56.  

\smallskip
\ref{[Sh1]} S. Shelah, Proper forcing, Lec.
Notes in Math. 940, Springer-Verlag (1982).

\smallskip
\ref{[Sh2]} S. Shelah, Cardinal Arithmetic,
Oxford Univ. Press (1994).

\smallskip
\ref{[Sh3]} S. Shelah, Appendix:  on stationary
sets (to ``Classification of nonelementary
classesII).  Abstract elementary classes") in
Proc. of the USA-Israel Conference on
Classification theory, J. Baldwin ed., Lec.
Notes in Math. 1292, Springer (1987), 483-485.

\smallskip
\ref{[Sh4]} S. Shelah, On Successors of Singular
Cardinals, in Logic Colloquium 78, M. Boffa, D.
van Dalen and K. McAloon (eds.), North Holland
(1979), 357-380.

\smallskip
\ref{[Sh5]} S. Shelah, Diamonds and
uniformizations, J. Symb. Logic 49 (1984),
1022-1033.

\smallskip
\ref{[Sh6]} S. Shelah, Iterated forcing and
normal ideals on $\ome_1$, Is. Jour. of Math. 60
(1987), 345-380.

\smallskip
\ref{[Sh7]} S.Shelah,More on stationry coding, in Around Classification
Theory of Models ,Lec.Notes in Math.1182,(1986),230-233

\smallskip
\ref{[Sh8]} S.Shelah,Non structure theory,Oxford Univ. Press ([Sh-h]) ,
to appear.

\smallskip
\ref{[Sh9]} S.Shelah,The Generalized Continuum Hypothesis Revisited
 ([Sh 460]),to appear

\smallskip
\ref{[So]} R. Solovay, Real valued measurable
cardinals, in Axiomatic Set Theory, Proceedings
of Sym. in Pure Math. XIII, vol. 1 (1971),
387-428.

\smallskip
\ref{[St]} J. Steel, The Core Model Iterability
Problem, to appear.

\smallskip
\ref{[St-V]} J. Steel and R. Van Wesep, Two
consequence of determinacy consistent with
choice, Trans. Am. Math. Soc.272 (1), (1982),
67-85.  

\smallskip
\ref{[Wo1]} H. Woodin, Some consistency results
in ZF using AD, Cabal Seminar, Lect. Notes in   
Math. 1019, Springer, 172-199.

\end